\documentclass[12pt,reqno]{amsart}
\usepackage{enumerate}
\usepackage{amsmath, amsfonts, amssymb, amsthm}
\def\Z{\Bbb Z}

\def\bg{\bigg}
\def\({\bg(}
\def\){\bg)}
\def\]{]\!]}
\def\[{[\![}

\def\f{\frac}
\def\mo{{\rm{mod}\ }}

\def\ls{\leqslant}
\def\gs{\geqslant}
\def\se {\subseteq}
\def\sm{\setminus}

\def\al{\alpha}

\def\eq{\equiv}

\def\Proof{\noindent{\it Proof}}
\theoremstyle{plain}
\newtheorem{theorem}{Theorem}
\newtheorem{lemma}{Lemma}

\newtheorem{conjecture}{Conjecture}
\theoremstyle{definition}

\theoremstyle{remark}

 \vspace{4mm}

\begin{document}
 \baselineskip=17pt
\hbox{Acta Arith., 140(2009), no.\,4, 329--334.}
\medskip

\title
[On Bialostocki's conjecture for zero-sum sequences] {On Bialostocki's conjecture for zero-sum sequences}

\author
[Song Guo and Zhi-Wei Sun] {Song Guo and Zhi-Wei Sun*}

\thanks{*Supported by the National Natural Science Foundation (grant 10871087) of China}

\address{(Song Guo) Department of Mathematics, Huaiyin Normal College,
Huaian 223300, People's Republic of China}
\email{guosong77@hytc.edu.cn}
\address {(Zhi-Wei Sun) Department of Mathematics, Nanjing
University, Nanjing 210093, People's Republic of China}
\email{zwsun@nju.edu.cn}


 \begin{abstract} Let $n$ be a positive even integer, and let $a_1, \ldots , a_n$ and $w_1,\ldots , w_n$ be
integers satisfying $\sum_{k=1}^na_k\eq\sum_{k=1}^nw_k
\eq 0\ (\mo n)$.  A conjecture of Bialostocki states that
there is a permutation $\sigma$ on $\{1,\ldots,n\}$ such that
$\sum_{k=1}^nw_ka_{\sigma(k)}\eq 0\ (\mo n)$. In this paper we confirm the conjecture when
$w_1,\ldots, w_n$ form an arithmetic progression with even
common difference.
\end{abstract}

\maketitle

\section{Introduction}

A finite sequence $S$ of terms from an (additive) abelian group is said to have zero-sum if the sum
of the terms of $S$ is zero. In 1961 P. Erd\H{o}s, A. Ginzburg and
A. Ziv \cite{EGZ} proved that any sequence of $2n - 1$ terms from an
abelian group of order $n$ contains an $n$-term zero-sum
subsequence. This celebrated EGZ theorem is is an important result in
combinatorial number theory and it has many different generalizations \cite{G1,H,S,T} including
Sun's recent extension involving covering systems.

The following theorem is called the weighted EGZ theorem. It was
conjectured by Y. Caro \cite{C} and proved by D. J. Grykiewicz
\cite{G}.

\begin{theorem}[Weighted EGZ Theorem] Let $n$ be a positive integer
and let $w_1,\ldots,w_n\in\Z_n=\Z/n\Z$ with
$\sum_{k=1}^n w_k = 0$. If $a_1, a_2, \ldots , a_{2n-1}$ is a
sequence of elements from $\Z_n$, then $\sum_{k=1}^nw_ka_{j_k}= 0$
for some distinct $j_1,\ldots,j_n\in\{1, \ldots 2n-1\}$.
\end{theorem}

Recently Bialostocki raised the following challenging conjecture.

\begin{conjecture}[{\rm Bialostocki \cite[Conjecture 14]{Bi}}]
Let $n$ be a positive even integer.
Suppose that $a_1,\ldots , a_n$ and $w_1, \ldots , w_n$ are
zero-sum sequences with terms from $\Z_n$.  Then there exists a
permutation $\sigma\in S_n$ such that $\sum_{k=1}^nw_ka_{\sigma(k)} = 0$, where
$S_n$ denotes the symmetric group of all permutations on $\{1,\ldots,n\}$.
\end{conjecture}

The conjecture has been verified for $n=2,4,6,8$. It fails for $n=3,5,7,\ldots$.
For example, $\{a_1, a_2, a_3\}=\{w_1, w_2,
w_3\}=\Z_3$ gives a counter-example for $n=3$.

In this paper we mainly establish the following result.

\begin{theorem} Let $n$ be a positive even integer,
and let $a_1, \ldots, a_n\in\Z$ with $\sum_{k=1}^na_k\eq 0\ (\mo n)$.
Then there exists a permutation $\sigma\in S_n$ such that
$\sum_{k=1}^nka_{\sigma(k)} \eq 0\ (\mo n/2)$.
Consequently, if $w_1,\ldots,w_n\in\Z$ form an arithmetic progression with even common difference,
then $\sum_{k=1}^nw_ka_{\sigma(k)} \eq 0\ (\mo n)$ for some $\sigma\in S_n$.
\end{theorem}

We are going to present two lemmas in the next section and then give our proof of Theorem 1.2 in Section 3.

\section{Two Lemmas}

\begin{lemma} Let $n=mq$ with $m,q\in\Z^+=\{1,2,3,\ldots\}$ and
$m\ge2$. Let $d\in\Z^+$ be a divisor of $q$, and let $a_1,\ldots,a_n\in\Z$.
Then, there is a partition $I_1,\cdots,I_m$ of $[1,n]=\{1,\ldots,n\}$
such that for each $s=1,\ldots,m$ we have $|I_s|=q$ and
$$d\ \big|\ \sum_{i\in I_s}a_i\ \Longrightarrow\ |\{a_i\ \mo d:\ i\in I_s\}|=1.$$
\end{lemma}
\Proof. By induction on $m$, it suffices to show that there exists an $I\se[1,n]$ with $|I|=q$
such that for each $J\in\{I,\,[1,n]\sm I\}$ we have $|\{a_j\ \mo d:\,j\in J\}|=1$ or
$\sum_{j\in J}a_j\not\eq0\ (\mo\ d)$. To achieve this we distinguish three cases.

{\it Case} 1. $|\{a_i\ \mo d:\,i\in[1,n]\}|=1$.

 In this case, $I=[1,q]$ works for our purpose.

{\it Case} 2. $|\{a_i\ \mo d:\,i\in[1,n]\}|=2$.

Suppose that
$$\{a_i\ \mo d:\, i\in[1,n]\}=\{r\mo d,\ r'\mo d\},$$
where $r,r'\in[0,d-1]$, $r\not\eq r'\ (\mo d)$, and
$a_i\eq r\ (\mo d)$ for at least $n/2$ values of $i\in[1,n]$.
Choose $I_0\se\{i\in[1,n]:\,a_i\eq r\ (\mo d)\}$ with $|I_0|=q\ls n/2$.
Let $i_0\in I_0$ and $j_0\in\bar I_0=[1,n]\sm I_0$ with $a_{j_0}\eq r'\ (\mo d)$.
When $\sum_{j\in\bar I_0}a_j\eq0\ (\mo d)$, we have both
$$\sum_{i\in (I_0\sm\{i_0\})\cup\{j_0\}}a_i\eq 0-r+r'\not\eq0\ (\mo d)$$
and
$$\sum_{j\in (\bar I_0\sm\{j_0\})\cup\{i_0\}}a_j\eq 0-r'+r\not\eq0\ (\mo d).$$
Thus, there always exists an $I\se[1,n]$ with $|I|=q$ such that
$$|\{a_i\ \mo d:\, i\in I\}|=1\ \mbox{or}\ \sum_{i\in I}a_i\not\eq0\ (\mo d),$$
and also $\sum_{j\in\bar I}a_j\not\eq0\ (\mo d)$.

{\it Case} 3. $|\{a_i\ \mo d:\,i\in[1,n]\}|>2$.

As $n\gs 2q\gs2q-1$, by the EGZ theorem there is an $I_0\se[1,n]$ with $|I_0|=q$ such that
$\sum_{i\in I_0}a_i\eq0\ (\mo q)$. For $\bar I_0=[1,n]\sm I_0$, we clearly have $|\bar I_0|=(m-1)q$.
Set $b=a_1+\cdots+a_n\eq\sum_{j\in \bar I_0}a_j\ (\mo\ q)$.

Suppose that  $a_j-a_i\eq 0$ or $b\ (\mo d)$ for any $i\in I_0$ and
$j\in \bar I_0$. Then
$$|\{a_i\ \mo d:\, i\in I_0)\}|\ls2\ \ \ \mbox{and}\ \ \ |\{a_j\ \mo d:\, j\in \bar I_0)\}|\ls2.$$
If $i_1,i_2\in I_0$, $j\in\bar I_0$ and $a_j\not\eq a_{i_1},a_{i_2}\
(\mo p)$, then $a_j-a_{i_1}\eq b\eq a_j-a_{i_2}\ (\mo d)$ and hence
$a_{i_1}\eq a_{i_2}\ (\mo d)$. So, if $|\{a_i\ \mo d:\,i\in I_0\}|=2$
then $\{a_j\ \mo d:\, i\in \bar I_0\}\se \{a_i\ \mo d:\,i\in I_0\}$
which contradicts $|\{a_i\mo d:\,i\in I_0\}|>2$. Similarly, if
$|\{a_j\mo d:\,j\in \bar I_0\}|=2$ then we also have a
contradiction. When $$|\{a_i\ \mo d:\, i\in I_0\}|=|\{a_j\mo d:\,j\in \bar I_0\}|=1,$$
we cannot have $|\{a_i\ \mo d:\,i\in[1,n]\}|>2$.

By the above, there are $i_0\in I_0$ and $j_0\in\bar I_0$ such that
$$a_{j_0}-a_{i_0}\not\eq 0,b\ (\mo d).$$ Set
$$I= (I_0\sm\{i_0\})\cup\{j_0\}\ \ \mbox{and}\ \ \bar I=[1,n]\sm I=(\bar I_0\sm\{j_0\})\cup\{i_0\}.$$
Then
$$\sum_{i\in I}a_i=\sum_{i\in I_0}a_i-a_{i_0}+a_{j_0}=0-a_{i_0}+a_{j_0}\not\eq0\ (\mo d)$$
and
$$\sum_{j\in \bar I}a_j=\sum_{j\in\bar I_0}a_j-a_{j_0}+a_{i_0}\eq b+a_{i_0}-a_{j_0}\not\eq0\ (\mo d).$$
Note that $|I|=q$ and $|\bar I|=(m-1)q$.

 Combining the above and using the induction argument, we see that
 the desired result holds for any $m=2,3,4,\ldots$.
\qed

\begin{lemma} Let $a_1,\ldots,a_n\in \Z$ with $n=p^{\al}$,
where $p$ is an odd prime and $\al$ is a positive integer. If
$\sum_{k=1}^n a_k \not\eq 0(\mo p)$ or $|\{a_k\ \mo\ p:\,k\in [1,n]\}|=1$, then
there exists a permutation $\sigma\in
S_n$ such that $\sum_{k=1}^nka_{\sigma(k)}\eq 0\ (\mo n)$.
\end{lemma}
\Proof. If $a:=\sum_{k=1}^n a_k \not\eq 0\ (\mo p)$, then
there is an $l\in[1,n]$ such that $al+\sum_{k=1}^n ka_k\eq0\ (\mo p)$
and hence
$$\sum_{k=1}^n ka_{\sigma(k)}\eq\sum_{k=1}^n(k+l)a_{k}\eq \sum_{k=1}^nka_k+la\eq 0\ (\mo p^{\al}),$$
where $\sigma(k)$ is the least positive residue of $k-l$ modulo $n$.

In the case $a_1\eq\cdots\eq a_n\ (\mo p)$, it is clear that
$$\sum_{k=1}^pka_k\eq a_1\sum_{k=1}^p k=a_1p\f{p+1}2\eq 0\ (\mo p).$$
Thus we have the desired result for $\al=1$.

Now let $\al>1$ and assume the desired result with $\al$ replaced by $\al-1$.
As mentioned above, the desired result holds if $\sum_{k=1}^n a_k \not\eq 0\ (\mo p)$.
Suppose that $a_1\eq\cdots\eq a_n\ (\mo p)$ and set $b_k=(a_k-a_1)/p$ for $k=1,\ldots,n$.
In light of Lemma 2.1, there exists a
partition $I_1\cup\cdots\cup I_p$ of $[1,n]$ with $|I_1|=\cdots=|I_p|=p^{\al-1}$ such that
for any $s=1,\ldots,p$ either $|\{b_k\ \mo p:\, k\in I_s\}|=1$ or $\sum_{k\in I_s}b_k\not\eq0\ (\mo\ p)$.
By the induction hypothesis,
there are one-to-one mappings $\sigma_s:[1,p^{\al-1}]\to I_s$
($s=1,\ldots,p$) such that
$$\sum_{k=1}^{p^{\al-1}}kb_{\sigma_s(k)}\eq 0 (\mo p^{\al-1})\ \ \mbox{for all}\
s=1\ldots,p.$$
For $s\in[1,p]$ and $t\in[1,p^{\al-1}]$ define $\sigma(p^{\al-1}(s-1)+t)=\sigma_s(t)$.
Then $\sigma\in S_n$ and
\begin{align*}
\sum_{k=1}^nka_{\sigma(k)}=&\sum_{k=1}^n ka_1+p\sum_{k=1}^kkb_{\sigma(k)}
\\=&\f{p^\al(p^\al+1)}2a_1+p\sum_{s=1}^p\sum_{t=1}^{p^{\al-1}}(p^{\al-1}(s-1)+t)b_{\sigma_s(t)}
\\\eq&p\sum_{s=1}^p\sum_{t=1}^{p^{\al-1}}tb_{\sigma_s(t)}\eq 0\ (\mo p^{\al}).
\end{align*}
This concludes the induction step and we are done. \qed

\section{Proof of Theorem 1.2}

\medskip
\noindent{\it Proof of Theorem 1.2}. We use induction on $\nu(n)$, the total number of prime divisors of $n$.

In the case $\nu(n)=1$, clearly $n=2$ and the desired result holds trivially.

Now let $\nu(n)>1$ and assume the desired result for those even positive integers
with less than $\nu(n)$ prime divisors.

{\it Case} 1. $n=2^{\al}$ for some $\al\gs2$.

By the EGZ theorem, there is an $I\se[1,n]$ with $|I|=n/2=2^{\al-1}$
such that $\sum_{i\in I}a_i\eq0\ (\mo\ 2^{\al-1})$. Note that for $\bar I=[1,n]\sm I$ we also have
$$\sum_{j\in\bar I}a_j=\sum_{k=1}^na_k-\sum_{i\in I}a_i\eq0\ (\mo 2^{\al-1}).$$
By the induction hypothesis,  for some one-to-one mappings $\sigma_0:[1,n/2]\to I$
 and $\sigma_1:[1,n/2]\to \bar I$ we have
$$2\sum_{k=1}^{2^{\al -1}}ka_{\sigma_0(k)}\eq 2\sum_{k=1}^{2^{\al -1}}ka_{\sigma_1(k)} \eq 0\ (\mo 2^{\al -1}).$$
Observe that
$$\sum_{k=1}^{2^{\al -1}}(2k-1)a_{\sigma_1(k)}
\eq 2\sum_{k=1}^{2^{\al -1}}ka_{\sigma_1(k)}-\sum_{j\in\bar I} a_j\eq 0\ (\mo 2^{\al -1}).$$
For $k\in[1,n/2]$ and $r\in[0,1]$ define $\sigma(2k-r)=\sigma_r(k)$. Then
$\sigma \in S_n$  and
$$\sum_{j=1}^{n}j a_{\sigma(j)}= 2\sum_{k=1}^{2^{\al-1}}ka_{\sigma_0(k)}
+\sum_{k=1}^{2^{\al-1}}(2k-1)a_{\sigma_1(k)}\eq 0\ (\mo  2^{\al-1}).$$
Thus we have the desired result for $n=2^\al$.

{\it Case} 2. $n$ has an odd prime divisor $p$.

Write $n=p^\al m$ with $\al,m>0$ and $p\nmid m$.
With the help of Lemma 2.1 there is a partition $I_1\cup\cdots I_m$ of $[1,n]$
with $|I_1|=\cdots=|I_m|=p^{\al}$ such that for each $s=1,\ldots,m$ either
$|\{a_i\ \mo p:\,i\in I_s\}|=1$ or $\sum_{i\in I_s}a_i\not\eq0\ (\mo p)$.
Combining this with Lemma 2.2, we see that for each $s\in[1,m]$ there is a one-to-one mapping
$\sigma_s:[1,p^\al]\to I_s$ such that $\sum_{t=1}^{p^\al}ta_{\sigma_s(t)}\eq0\ (\mo\ p^\al)$.

Set $b_s=\sum_{k\in I_s}a_k$ for $s=1,\ldots,m$. Then
$$\sum_{s=1}^mb_s=\sum_{k\in I_1\cup\cdots\cup I_m}a_k=\sum_{k=1}^n a_k\eq0\ (\mo m).$$
As $2\mid m$ and $\nu(m)<\nu(n)$, by the induction hypothesis, for some $\tau\in S_m$ we have
$$2\sum_{s=1}^msb_{\tau(s)}\eq 0\ (\mo m)$$
and hence
$$2\sum_{s=1}^m\sum_{t=1}^{p^\al}sa_{\sigma_{\tau(s)}(t)}
=2\sum_{s=1}^msb_{\tau(s)}\eq0\ (\mo m).$$
Note also that
$$\sum_{s=1}^m\sum_{t=1}^{p^\al}ta_{\sigma_{\tau(s)}(t)}
=\sum_{s=1}^m\sum_{t=1}^{p^\al}ta_{\sigma_s(t)}\eq0\ (\mo p^\al).$$
Therefore
$$2\sum_{s=1}^m\sum_{t=1}^{p^\al}(p^\al s+mt)a_{\sigma_{\tau(s)}(t)}\eq0\ (\mo p^\al m).$$
As $p^\al$ is relatively prime to $m$,
$$\{p^\al s+mt:\, s\in[1,m]\ \mbox{and}\ t\in[1,p^\al]\}$$
is a complete system of residues modulo $n=p^\al m$.
For any $k\in[1,n]$, there are unique $s\in[1,m]$ and $t\in[1,p^\al]$ such that
$k\eq p^\al s+mt\ (\mo\ n)$, and we define $\sigma(k)=\sigma_{\tau(s)}(t)$.
Then $\sigma\in S_n$ and also
$$2\sum_{k=1}^nka_{\sigma(k)}\eq 0 (\mo n).$$
This concludes the induction step.

 In view of the above, we have completed the proof of Theorem 1.2. \qed


\begin{thebibliography}{99}

\bibitem {Bi} A. Bialostocki, \textit{Some problems in view of
recent developments of the Erd\H{o}s-Ginzburg-Ziv Theorem},
Integers {\bf 7} (2007), no.2, {\#} A07, 10 pp (electronic).

\bibitem {C} Y. Caro,
\textit{Zero-sum problems--a survey}, Discrete Math. {\bf 152}
(1996), 93--113.

\bibitem {EGZ} P. Erd\H os, A. Ginzburg and A. Ziv,
\textit{Theorem in additive number theory}, Bull. Res. Council
Israel {\bf 10F} (1961), 41--43.

\bibitem {G} D. J. Grynkiewicz, \textit{A weighted Erd\H{o}s-Ginzburg-Ziv Theorem},
Combinatorica {\bf 26} (2006), 445--453.

\bibitem {G1} D. J. Grynkiewicz, \textit{On the number of m-term zero-sum
subsequences}, Acta Arith. {\bf 121} (2006), 275--298.

\bibitem {H} Y. O. Hamidoune, O. Ordaz and A. Ortu\~{n}o, \textit{On a combinatorial
theorem of Erd\H{o}s, Ginzburg and Ziv}, Combin. Probab. Comput.
{\bf 7} (1998), 403--412.

\bibitem {S} Z. W. Sun, \textit{Zero-sum problems for abelian $p$-groups and covers of the integers by residue classes},
Israel J. Math. {\bf 170} (2009), 235--252.

\bibitem {T} R. Thangadurai, \textit{Non-canonical extensions of Erd\H{o}s-Ginzburg-Ziv
theorem}, Integers {\bf 2} (2002), \#A07, 14 pp (electronic).


\end{thebibliography}
\end{document}